\documentclass[10pt]{amsart} 

\usepackage{preamble} 

\usepackage{fancyhdr}
\title{Effective bounds on characterising slopes for all knots} 
\author{Patricia Sorya}
\author{Laura Wakelin}
\thanks{The first author was supported by FRQNT doctoral research grant [305903]; the second author was partially funded by the EPSRC Centre for Doctoral Training in Geometry and Number Theory at the Interface: LSGNT [EP/S021590/1]. } 
\date{} 

\begin{document} 

\begin{abstract}
\vspace{-1em}
    A slope $p/q$ is characterising for a knot $K \subset \mathbb{S}^3$ if the orientation-preserving homeomorphism type of the manifold $\mathbb{S}^3_K(p/q)$ obtained by performing Dehn surgery of slope $p/q$ along $K$ uniquely determines the knot $K$. We combine new applications of results from hyperbolic geometry with previous individual work of the authors to determine, for any given knot $K$, an explicit bound $\mathcal{C}(K)$ such that $|q| > \mathcal{C}(K)$ implies that $p/q$ is a characterising slope for $K$. Furthermore, we find an optimal such $\mathcal{C}(K)$ for certain satellite knots with winding number zero patterns. 
\end{abstract} 

\maketitle 

\vspace{-1em}

\section{Introduction} 
\label{sec:introduction}

Given a knot $K \subset \mathbb{S}^3$, one can perform \emph{Dehn surgery} of slope $p/q \in \mathbb{Q} \cup \{1/0\}$ on $K$ to produce a new 3-manifold $\mathbb{S}^3_K(p/q) = \mathbb{S}^3_K \cup \nu(K)$, where $\mathbb{S}^3_K := \mathbb{S}^3 \setminus \interior(\nu(K))$ is the exterior of $K$ and the fraction $p/q$ specifies how to glue back in the solid torus $\nu(K) \cong \mathbb{S}^1 \times \mathbb{D}^2$. 
When $p/q=1/0$, it is impossible to uniquely determine the knot $K$ from the oriented homeomorphism type of the manifold $\mathbb{S}^3_K(1/0) \cong \mathbb{S}^3$.
However, this is not true for a general slope $p/q \in \mathbb{Q}$. 
We say that a slope $p/q \in \mathbb{Q}$ is \emph{characterising} for a knot $K \subset \mathbb{S}^3$ if the existence of an orientation-preserving homeomorphism $\mathbb{S}^3_K(p/q) \cong \mathbb{S}^3_{K'}(p/q)$ for $K' \subset \mathbb{S}^3$ implies that $K=K'$. 

The Dehn surgery characterisation problem asks which slopes $p/q \in \mathbb{Q}$ are characterising for a given knot $K$. 
So far, it has been solved for the unknot, the trefoil knots and the figure-eight knot, for which all slopes in $\mathbb{Q}$ are characterising \cite{kmos, os}, as well as for an infinite family of composite knots, for which the set of characterising slopes is precisely $\mathbb{Q}\setminus\mathbb{Z}$ \cite{sorya}. 
An advance towards answering this question for knots in general is a result of the first author, building on prior work on the subject \cite{lackenby, mccoy}, which says that for any knot $K \subset \mathbb{S}^3$, there exists a constant $\mathcal{C}(K)$ such that every slope $p/q$ with $|q| > \mathcal{C}(K)$ is characterising for $K$ \cite{sorya}. 
Whilst the proof of existence of $\mathcal{C}(K)$ is non-constructive in general, the present paper provides an explicit value of $\mathcal{C}(K)$ that depends only on the JSJ decomposition of the exterior of $K$.

\begin{theorem}
\label{theorem:main} 
    Let $K \subset \sphere$ be a knot with exterior $\sphere_K$. 
    Then the JSJ decomposition of $\sphere_K$ -- namely, the geometry of the JSJ pieces of $\sphere_K$, together with the gluing maps between them -- explicitly determines a constant $\mathcal{C}(K)$ such that if $|q| > \mathcal{C}(K)$, then $p/q$ is a characterising slope for $K$.
\end{theorem}

For torus knots $T_{a,b}$, the constant $\mathcal{C}(K)$ can be realised as $\max\{8,|a|,|b|\}$ \cite{mccoy}. 
For hyperbolic knots, the second author showed that $\mathcal{C}(K)$ can be constructed using the systole of the hyperbolic knot exterior \cite{wakelin}. 
For satellite knots, explicit bounds are known for prime knots whose exteriors consist only of Seifert fibred JSJ pieces and for all composite knots \cite{sorya}; for the latter, we can take $\mathcal{C}(K)=1$. 

This leaves the case where $K$ is a prime satellite knot whose exterior contains at least one hyperbolic JSJ piece. 
In this paper, we construct a value for $\mathcal{C}(K)$ for such knots which only depends on the geometry of the hyperbolic JSJ pieces of $\sphere_K$ and how they are glued within the full JSJ decomposition.

\subsection{Determining a value for \texorpdfstring{$\mathcal{C}(K)$}{C(K)}}

The remaining case of $\mathcal{C}(K)$ is realised by taking the maximum of three geometric constants, $Q(K)$, $R(K)$ and $S(K)$, which we now introduce and define. 

\begin{theorem}
\label{theorem:C(K)}
    Let $K$ be a prime satellite knot whose exterior is not a graph manifold.
    
    If $|q|> \max\{Q(K), R(K), S(K)\}$, then $p/q$ is a characterising slope for $K$.
\end{theorem}

To define these constants, consider the JSJ decomposition of $\sphere_K$. 
Each JSJ piece of $\sphere_K$ is homeomorphic to the exterior of a unique link $L = L_0 \cup L_1 \cup \ldots \cup L_{m-1}$ with a distinguished component $L_0$ such that $L \setminus L_0 = U^{m-1}$ is a (possibly empty) unlink \cite[Proposition~2.4]{budney}. 
Denote by $\mu_i$ the meridian of the $i^\text{th}$ component of $U^{m-1}$ for $i = 1, \ldots, m-1$. 
When $L$ is hyperbolic, we set
\begin{align*}
    \mathfrak{q}(L) &:= 
    \Bigg\lfloor \sqrt{6\sqrt{3} \bigg( 1.9793 \frac{2\pi}{\sys(\sphere_L)} + 28.78 \bigg)} \ \Bigg\rfloor, \\
    \mathfrak{r}(L) &:= 
    \Bigg\lfloor\sqrt{3}\max\{0, l(\mu_i) \;|\; L_i \subset L, i\neq 0\}\Bigg\rfloor, \\
    \mathfrak{s}(L) &:=  
    \Bigg\lfloor\sqrt{6\sqrt{3}\bigg(\frac{2\pi}{\sys(\sphere_L)}+28.78\bigg)}\Bigg\rfloor,
\end{align*}
where $\sys(\sphere_L)$ denotes the length of the shortest geodesic in $\sphere_L$ and $l(\mu_i)$ is the length of $\mu_i\subset\partial\sphere_L$; these quantities are described in more detail in Subsection \ref{subsection:hyperbolic}.

Let $\mathcal{X}$ denote the set of hyperbolic JSJ pieces of $\mathbb{S}^3_K$. 
For each $X \in \mathcal{X}$, let $L_X$ denote the unique hyperbolic link corresponding to $X$ in the satellite construction of $K$.

To define $Q(K)$, we only need to consider one JSJ piece of $\sphere_K$. 
If $K$ is a cable of a knot $\widehat{K}$, let $Y$ be the JSJ piece of $\sphere_{\widehat{K}}$ containing $\partial \sphere_{\widehat{K}}$; otherwise, let $Y$ be the JSJ piece of $\sphere_K$ containing $\partial\sphere_K$. 

\begin{definition}
    Define $Q(K) := \max\{34, \mathfrak{q}(L_Y)\}$ if $Y \in \mathcal{X}$ and set $Q(K) = 0$ otherwise. 
\end{definition} 

For the constant $R(K)$, we use $\mathfrak{r}(L_X)$ for $X \in \mathcal{X}$, which can only contribute non-trivially if $|\partial X|\geq2$. 

\begin{definition}
    Define $R(K) := \max\{1, \mathfrak{r}(L_X) \ | \ X \in \mathcal{X}\}$. 
\end{definition} 

Finally, the constant $S(K)$ involves every $X \in \mathcal{X}$ which does not contain $\partial \sphere_K$. 

\begin{definition}
    Define $S(K) := \max\{32, \mathfrak{s}(L_X) \ | \ X \in \mathcal{X}, \ \partial \sphere_K \not \subset X\}$.
\end{definition} 

Combining these new constants with the ones from previous work, we obtain the following explicit description of a constant $\mathcal{C}(K)$ that realises Theorem \ref{theorem:main} for all knots.

\begin{theorem}
\label{theorem:allcases} 
    Let $K \subset \sphere$ be a knot with exterior $\sphere_K$.
    \begin{enumerate}[label={(\roman*)}]
    \item If $K$ is the unknot, set $\mathcal{C}(K)=0$. 
        \item If $K$ is a composite knot, set $\mathcal{C}(K)=1$. 
        \item If $K$ is a prime knot and $\sphere_K$ is a graph manifold, express $K$ as an iterated cable $C_{r_n, s_n}\ldots C_{r_2, s_2}C_{r_1, s_1}(J)$, where $s_i$ is the winding number of $C_{r_i, s_i}$ and $J$ is either a torus knot or a composite knot, and set 
        \[\mathcal{C}(K) = 
        \begin{cases}
            \max\{8, |a|, |b|\}  &\text{ if $J$ is a torus knot $T_{a,b}$ and }n=0,\\
            \max\{8,|s_1|,|r_1|+|a|, |r_1|+|b|\} &\text{ if $J$ is a torus knot $T_{a,b}$ and }n=1,\\
            \max\{|r_1|+|a|, |r_1|+|b|\} &\text{ if $J$ is a torus knot $T_{a,b}$ and }n \geq 2,\\
            2 &\text{ if $J$ is a composite knot}.\\
        \end{cases}
        \]      
    \item If $K$ is a prime knot and $\sphere_K$ is not a graph manifold, set \[\mathcal{C}(K)=\max\{Q(K), R(K), S(K)\}.\]
    \end{enumerate}
    If $|q| > \mathcal{C}(K)$, then $p/q$ is a characterising slope for $K$. 
\end{theorem}

\begin{proof}
   The first case is due to \cite{kmos}; the second case is \cite[Theorem 2]{sorya}; the third case comes from \cite[Theorem 1.1]{mccoy} and \cite[Theorem 7.1]{sorya}; the fourth case is \cite[Theorem 1.3]{wakelin} together with Theorem \ref{theorem:C(K)}.  
\end{proof} 

The most challenging part of the proof of Theorem \ref{theorem:C(K)} lies in showing that the denominator bound \mbox{$|q|>\max\{R(K),S(K)\}$} ensures that any orientation-preserving homeomorphism \mbox{$f: \mathbb{S}^3_K(p/q) \to \mathbb{S}^3_{K'}(p/q)$} restricts to one between the \emph{surgered pieces} -- the JSJ pieces containing the surgery curves -- thereby realising \cite[Proposition~4.1]{sorya}. 
This is encompassed in the following result. 

\begin{proposition}
\label{proposition:swapping}
    Let $K$ be a prime satellite knot whose exterior is not a graph manifold. Suppose that there is an orientation-preserving homeomorphism $f: \mathbb{S}^3_K(p/q) \to \mathbb{S}^3_{K'}(p/q)$ for some knot $K'$.
    
    If $|q| > \max\{R(K), S(K)\}$, then $f$ maps the surgered piece of $\mathbb{S}^3_K(p/q)$ to the surgered piece of $\mathbb{S}^3_{K'}(p/q)$.
\end{proposition}

Once we have proved Proposition \ref{proposition:swapping}, Theorem \ref{theorem:C(K)} will then follow from the techniques of \cite{sorya} when the surgered piece is Seifert fibred and from the bound $|q| > Q(K)$ previously introduced in \cite{wakelin} when the surgered piece is hyperbolic.

\subsection{Winding number zero}

Note that in general, the bound provided by Theorem \ref{theorem:C(K)} may not be optimal. 
However, in certain cases of prime satellite knots whose exteriors are not graph manifolds, we can use a slightly different approach via a new constant $T(K)$ to find a refined value for $\mathcal{C}(K)$.

\begin{theorem}
\label{theorem:C(K)usingt}
    Let $K$ be a satellite knot such that for every choice of satellite description $K=P(J)$, the pattern $P$ has winding number zero. 
    
    If $|q|> \max\{Q(K),T(K)\}$, then $p/q$ is a characterising slope for $K$. 
\end{theorem} 

With this strategy, we may be able to improve the bound used to ensure that an orientation-preserving homeomorphism \mbox{$f: \mathbb{S}^3_K(p/q) \to \mathbb{S}^3_{K'}(p/q)$} restricts to one between surgered pieces. Namely, instead of assuming that $|q|>\max\{R(K),S(K)\}$, we write $K=P(J)$ and use the fact that the pattern $P$ has winding number zero to deduce that $f$ restricts to a homeomorphism between the surgered pieces unless $|p|=1$. If we write $K'=P'(J')$ and suppose that $f(\mathbb{S}^3_{J}) = \mathbb{S}^3_{P'}(\pm1/q)$ and $f(\mathbb{S}^3_P(\pm1/q)) = \mathbb{S}^3_{J'}$, we can express these $\pm1/q$-surgeries as Rolfsen $\mp q$-twists. We obstruct this by introducing a constant $T(K)$ encoding a certain unknotting property of the companion $J$. This leads to the following proposition. 

\begin{proposition}
\label{proposition:swappingusingt}
    Let $K$ be a satellite knot such that for every choice of satellite description $K=P(J)$, the pattern $P$ has winding number zero. Suppose that there is an orientation-preserving homeomorphism $f: \mathbb{S}^3_K(p/q) \to \mathbb{S}^3_{K'}(p/q)$ for some knot $K'$.
    
    If $|q| > \max\{2, T(K)\}$, then $f$ maps the surgered piece of $\mathbb{S}^3_K(p/q)$ to the surgered piece of $\mathbb{S}^3_{K'}(p/q)$.
\end{proposition}

Combining Proposition \ref{proposition:swappingusingt} with the second author's previous work \cite{wakelin} leads to Theorem \ref{theorem:C(K)usingt}. Furthermore, we show that sometimes this bound is optimal. 

\subsection{Outline} 
The paper is organised as follows. 
Section \ref{sec:preliminaries} reviews the structure of satellite knots, as well as JSJ decompositions of knot exteriors and their Dehn fillings, and recaps useful results from hyperbolic geometry. 
Section \ref{sec:realising} presents the proof of Theorem \ref{theorem:C(K)}. 
Section \ref{sec:windingzero} contains the proof of Theorem \ref{theorem:C(K)usingt}. 
Section \ref{sec:applications} comprises applications of our main results: we give examples showing how to compute the bound from Theorem \ref{theorem:C(K)} and investigate some special cases where the bound from Theorem \ref{theorem:C(K)usingt} is optimal.

\subsection{Acknowledgements} 
We are grateful to Imperial College London, the Max-Planck-Institut f\"ur Mathematik (MPIM) and the Université du Québec à Montréal (UQAM) for their generous support during this joint work, as well as to the University of Oxford and the Rényi Alfréd Matematikai Kutatóintézet for providing us with the opportunity to collaborate in stimulating environments. 
We thank the anonymous referee for their thorough review and invaluable comments. 
We extend our thanks to Steven Boyer, Marc Lackenby, Duncan McCoy and Steven Sivek for their continued advice and support throughout this project. 
We would also like to acknowledge Vitalijs Brejevs and Diego Santoro for interesting discussions, as well as Nathan Dunfield for advice on an earlier version of this work and Dan Radulescu for coding suggestions. 
\section{Preliminaries}
\label{sec:preliminaries} 

In this section, we review satellite knots, JSJ decompositions and hyperbolic geometry in the context of this article.

\subsection{Satellite knots} 

Let $L = L_0 \cup L_1 \cup \ldots \cup L_{m-1}$ be an $m$-component link in $\sphere$, where each component $L_i$ is a knot. 
We denote by $\sphere_L$ the exterior of $L$ in $\sphere$: the manifold obtained by removing the interior of a closed tubular neighbourhood $\nu(L)$ of $L$ in $\sphere$. 

A \emph{slope} on a boundary component $\partial\nu(L_i)$ of $\sphere_L$ is an isotopy class of essential simple closed curves representing an element of $H_1(\partial\nu(L_i); \mathbb{Z})$ up to sign. 
The \emph{meridian} $\mu_i$ on $\partial\nu(L_i)$ is the unique slope that bounds a disc in $\nu(L_i)$. 
Let $M_i$ be the manifold obtained by Dehn filling, for each $j \neq i$, the boundary component $\partial\nu(L_j)$ of $\sphere_L$ along $\mu_j$, so that $M_i \cong \sphere_{L_i}$.
The \emph{longitude} $\lambda_i$ on $\partial\nu(L_i) = \partial M_i$ is the unique slope that is trivial in $H_1(M_i; \mathbb{Z})$. 
Fixing $\{\mu_i, \lambda_i\}$ as a basis for $H_1(\partial\nu(L_i); \mathbb{Z})$, the \emph{slope} $p/q\in\mathbb{Q} \cup \{1/0\}$ refers to $p\mu_i+q\lambda_i$ up to sign, where $p, q$ are coprime. 

We write $P(J)$ for the satellite knot with pattern $P$ and non-trivial companion knot $J$. 
The pattern $P$ can be described as a 2-component link $Q \cup U$, where $U$ is the unknot and $Q$ is a knot inside the solid torus $\sphere_U$ which is neither its core nor a local knot. 
The  satellite knot $P(J)$ is obtained by \emph{splicing} $U$ with $J$; the exterior of $P(J)$ is the \emph{splice} of $\sphere_P$ and $\sphere_J$ along $\partial \nu (U)$ and $\partial \nu (J)$, i.e. the meridian $\mu_U$ and longitude $\lambda_U$ of $U$ are identified with the longitude $\lambda_J$ and meridian $\mu_J$ of $J$, respectively. 

The \emph{winding number} of the pattern $P = Q \cup U$ is the absolute value of the algebraic intersection number between $Q$ and an essential disc in the solid torus $\sphere_U$. 
If $Q$ happens to be unknotted, then we say that $P$ is an \textit{unknotted pattern}. In this case, we can swap the components of $P$ via an isotopy and observe that the definition of the winding number respects this exchange for homological reasons. 
If a pattern $P = Q \cup U$ is such that $Q$ intersects an essential disc in $\sphere_U$ geometrically once, then $P$ is called a \emph{composing pattern} and $P(J)$ is precisely the \emph{composite knot} $Q \# J$ (also called the \emph{connected sum} of $Q$ and $J$). 
If a pattern $P = Q \cup U$ is such that $Q$ is isotopic in $\mathbb{S}^3_U$ to a torus knot, then $P$ is a \emph{cable pattern} or \emph{cable link} and $P(J)$ is a \emph{cable knot}.

\subsection{JSJ decompositions}

Recall that for any compact orientable irreducible 3-manifold $M$ whose boundary is a (possibly empty) union of tori, there is a minimal collection $\mathbf{T}$ of properly embedded disjoint essential tori such that each component of $M \setminus \mathbf{T}$ is either a hyperbolic 3-manifold or a Seifert fibre space; such a collection is unique up to isotopy \cite{js, j}. 
The \emph{JSJ decomposition} of $M$ is
\[M = M_0 \cup M_1 \cup \ldots \cup M_k,\]
where each $M_i$ is the closure of a component of $M \setminus \mathbf{T}$ and the gluing data is specified by self-homeomorphisms of each torus in $\mathbf{T}$. 
Each $M_i$ is called a \emph{JSJ piece} of $M$ and each torus in the collection $\mathbf{T}$ is called a \emph{JSJ torus} of $M$. 
If all of the JSJ pieces of $M$ are Seifert fibred, then $M$ is said to be a \emph{graph manifold}. 
For any homeomorphism $f: M \rightarrow N$ between compact orientable irreducible 3-manifolds, we may assume that, possibly after a homotopy, $f$ restricts to a homeomorphism between JSJ pieces. 
If $T$ is a JSJ torus of $M$ with gluing data given by a homeomorphism $\phi: T \rightarrow T$, then the gluing data along $f(T) \subset N$ is given by $f|_T \circ \phi \circ f^{-1}|_{f(T)}$. 
If the homeomorphism $f$ is orientation-preserving, then it also preserves induced orientations when restricted to JSJ pieces and JSJ tori. 

We will be most interested in the case where $M$ is the exterior $\mathbb{S}^3_K$ of a knot $K$. 
The unique JSJ piece of $\sphere_K$ containing the boundary of $\nu (K)$ is said to be the \emph{outermost} JSJ piece of $\sphere_K$. 

When $K$ is a torus knot or a hyperbolic knot, its exterior $\mathbb{S}^3_K$ contains no JSJ tori. 
When $K$ is a satellite knot, each JSJ torus splits $\sphere_K$ into a pattern space $\sphere_P$ and a knot exterior $\sphere_J$ such that $K$ can be described as $P(J)$ \cite[Lemma 2.1]{sorya}. 
In particular, this description may not be unique. 

The JSJ pieces of a knot exterior take on one of four special types, which are all homeomorphic to the exterior of a certain type of link in $\sphere$. 
The isotopy class of this link is uniquely determined by gluing maps arising in the satellite construction. 

\begin{theorem}[{\cite[Proposition 2.4, Theorem 4.18]{budney}}]
\label{theorem:budney}
    Let $K$ be a non-trivial knot and let $Y$ be a JSJ piece of the exterior of $K$. 
    Then $Y$ is homeomorphic to the exterior of a link $L = L_0 \cup L_1 \cup \ldots \cup L_{m-1}$ with a distinguished component $L_0$ such that $L \setminus L_0 = U^{m-1}$ is the unlink.
    
    Furthermore, this link $L$ is unique up to isotopy, given the condition that for each $i \neq 0$, the gluing map used in the satellite construction of $K$ corresponds to splicing $L_i$ with a non-trivial knot $J_i$.
    
    This unique link $L$ corresponding to $Y$ can be classified into one of the following four types.
    \begin{enumerate}[label={(\roman*)}]
        \item $L$ is a torus knot $T_{a,b}$, 
        i.e. $\sphere_L$ is a Seifert fibre space whose base orbifold is a disc with two cone points of orders $|a|$ and $|b|$. 
        \item $L$ is a cable link $C_{r,s}$,
        i.e. $\sphere_L$ is a Seifert fibre space whose base orbifold is an annulus with one cone point of order $|s|$ (where $s$ is the winding number of the cable pattern). 
        \item $L$ is a composing link, 
        i.e. $\sphere_L$ is a Seifert fibre space with at least three boundary components whose base orbifold is a planar surface with no cone points. 
        \item $L$ is a hyperbolic link, 
        i.e. $\sphere_L$ is a hyperbolic 3-manifold.
    \end{enumerate}
\end{theorem} 

This is summarised in the table below. 

\begin{center}
    \begin{table}[htbp!]
    \renewcommand{\arraystretch}{1.5}
        \begin{tabular}{|r|p{2.75cm}|p{4.1cm}|p{2cm}|p{4.1cm}|} \hline 
            Case & Link $L$ & JSJ piece $Y=\mathbb{S}^3_L$ & SFS orbifold & $K$ when $Y$ is outermost \\ \hline \hline 
            (i) & $T_{a,b}$ torus knot & $\mathbb{S}^3_{T_{a,b}}$ torus knot exterior & $\mathbb{D}^2(|a|,|b|)$ & $T_{a,b}$ torus knot \\ 
            (ii) & $C_{r,s}$ cable link & $\mathbb{S}^3_{C_{r,s}}$ cable space & $\mathbb{A}^2(|s|)$ & $C_{r,s}(\widehat{K})$ cable knot \\ 
            (iii) & $L$ composing link & $\mathbb{S}^3_L$ composing space & $\Sigma, |\partial \Sigma| \geq 3$ & $K_1 \# K_2$ composite knot \\ 
            (iv) & $L$ hyperbolic link & \mbox{$\mathbb{S}^3_L$ hyperbolic link exterior} & -- & \mbox{$K$ knot of hyperbolic type} \\ \hline 
        \end{tabular}
        \caption{JSJ pieces of the exterior of a knot $K$. }
        \label{table:budney}
    \end{table} 
\end{center}

\vspace{-2em}

The distinguished component $L_0$ of the link $L$ in the statement of Theorem \ref{theorem:budney} is said to be \emph{outermost}; similarly, the \emph{outermost} boundary component of $Y = \sphere_L$ refers to $\partial \nu (L_0) \subset \partial Y$.

\subsection{Dehn surgery}

Let $M$ be a 3-manifold and let $T_{0}, \ldots, T_{m-1}$ denote the toroidal boundary components of $\partial M$. 
For fixed bases $\{\mu_{i}, \lambda_{i}\}$ for each $H_1(T_{i}; \Z), i=0, \ldots, m-1$, let
\[M(T_{i}; p_{i}/q_{i})\]
denote the manifold obtained by Dehn filling $M$ along a simple closed curve representing $p_{i} \mu_{i}+q_{i} \lambda_{i}$ up to sign on $T_{i}$. 
If it is clear from context which boundary component of $M$ is filled, then we may simply write $M(p/q)$.

If $M$ is a link exterior $\sphere_L$ and $L_{0}, \ldots, L_{m-1}$ are the components of $L$, we may write \[\sphere_L(L_{i} ; p_{i}/q_{i})\]
instead of $\sphere_L(\partial\nu (L_{i}); p_{i}/q_{i})$. Furthermore, the slopes $p_{i}/q_{i}$ are assumed to be expressed with respect to the basis of $H_1(\partial \nu (L_{i}); \Z)$ given by the meridian and longitude of $L_{i}$, as defined earlier.

Performing Dehn surgery on a knot $K$ to obtain the manifold $\sphere_K(p/q)$ corresponds to Dehn filling the outermost JSJ piece of $\sphere_K$ along a slope $p/q$ on the boundary component corresponding to $\partial \nu (K)$.
The following proposition shows that if $|q|>2$, then the core of the surgery solid torus is contained inside a single JSJ piece of $\sphere_K(p/q)$. 
We call this piece the \emph{surgered piece}.

\begin{proposition}{\cite[Proposition~3.6]{sorya}}
\label{prop:JSJsurgery}
    Let $K$ be a non-trivial knot with exterior $\sphere_K$ and let \mbox{$Y_0 \cup Y_1 \cup \ldots \cup Y_k$} be the JSJ decomposition of $\sphere_K$, where $Y_0$ is the outermost piece of $\sphere_K$. 
    
    If $|q| > 2$, then the JSJ decomposition of $\sphere_K(p/q)$ takes one of the following forms: 
    \begin{enumerate}[label={(\roman*)}]
        \item $Y_0(\partial \nu(K); p/q) \cup Y_1 \cup Y_2 \cup \ldots \cup Y_k$,
        \item $Y_1(Y_0 \cap Y_1; p/qs^2) \cup Y_2 \cup \ldots \cup Y_k$,
    \end{enumerate} 
    where case (ii) occurs precisely when $K=C_{r,s}(\widehat{K})$ is a cable knot, $Y_1$ is the outermost piece of $\sphere_{\widehat{K}}$, $|s|\geq2$ and $|p-qrs|=1$.
\end{proposition}

Therefore the surgered piece can be described as $Y(p/qt^2)$, where either $Y=Y_0$ and $t=1$ or $Y=Y_1$ and $t > 1$.

\subsection{Hyperbolic geometry}
\label{subsection:hyperbolic}

When considering the hyperbolic JSJ pieces contributing to the constants in Theorem \ref{theorem:C(K)}, we will require some quantitative results from hyperbolic geometry. 

Recall that a hyperbolic $3$-manifold $M$ is one whose interior admits a complete finite-volume hyperbolic metric. A toroidal boundary component $T \subset \partial M$ is the boundary at infinity of a \emph{cusp} of the interior of $M$. Each of these has a well-defined \emph{maximal horocusp neighbourhood} $N(T)$, whose boundary $\partial N(T)$ inherits a unique Euclidean metric from the hyperbolic metric (see for instance \cite[Section 2]{lackenby}). 

\begin{definition}
    Let $\sigma$ be a slope on a toroidal boundary component $T \subset \partial M$.
    \begin{itemize}
        \item The \emph{area} $A(T)$ of $T$ is the Euclidean area of $\partial N(T)$.
        \item The \emph{length} $l(\sigma)$ of $\sigma$ is the Euclidean length of a geodesic representative of $\sigma$ on $\partial N(T)$. 
        \item The \emph{normalised length} $\hat{l}(\sigma)$ of $\sigma$ is given by $\hat{l}(\sigma) = l(\sigma)/\sqrt{A(T)}$. 
    \end{itemize}
\end{definition} 

Slope length is related to the geometry of Dehn fillings. 
The 6-theorem \cite{agol, lackenby_6} states that filling a hyperbolic 3-manifold along any slope $\sigma$ of length $l(\sigma)>6$ must give a hyperbolic 3-manifold. 

The lengths (or normalised lengths) of a pair of slopes on the same boundary component $T \subset \partial M$ can be related to the distance $\Delta$ between them (the absolute value of their algebraic intersection number), as presented in \cite[Lemma 4.2]{wakelin}. 

\begin{lemma}
\label{lemma:length_inequalities}
    Let $\sigma$ and $\sigma'$ be slopes on the same toroidal boundary component $T$ of a hyperbolic 3-manifold. 
    \begin{itemize} 
        \item The area $A(T)$ of $T$ satisfies the universal bound $A(T) \geq 2\sqrt{3}$. 
        \item The length $l(\gamma)$ of $\gamma$ satisfies $l(\gamma)l(\mu) \geq \Delta(\gamma, \mu) \cdot A(T)$. 
        \item The normalised length $\hat{l}(\gamma)$ of $\gamma$ satisfies $\hat{l}(\gamma)\hat{l}(\mu) \geq \Delta(\gamma, \mu)$. 
    \end{itemize} 
\end{lemma} 

\begin{proof}
    The first inequality comes from $\sqrt{3}$ being a universal lower bound for the volume of a maximal horocusp neighbourhood $N$ \cite[Theorem 1.2]{gabai} and the fact that the area of $\partial N$ is equal to twice the volume of $N$ \cite[Section 1]{gabai}. The others follow from \cite[Lemma 2.13]{cooper-lackenby} or \cite[Theorem 8.1]{agol}.
\end{proof}

Now we will move on to consider closed curves in the interior of $M$. 
We denote the length of any geodesic $\gamma \subset M$ with respect to the hyperbolic metric by $l(\gamma)$. 

\begin{definition}
    The \emph{systole} $\sys(M)$ of a hyperbolic 3-manifold $M$ is the length $l(\gamma_0)$ of a shortest geodesic $\gamma_0 \subset M$. 
\end{definition}

There is a quantitative relationship between the length $l(\gamma)$ of a simple geodesic $\gamma \subset M$ and the normalised length $\hat{l}(\mu)$ of its meridian $\mu$ considered as a slope on the boundary torus of the maximal horocusp neighbourhood of $M \setminus \nu (\gamma)$ corresponding to $\gamma$. 

\begin{theorem}{\cite[Corollary~6.13]{fps}}
\label{theorem:fps}
    Let $M$ be a hyperbolic 3-manifold and let $\gamma \subset M$ be a simple geodesic with meridian $\mu$.   
    If $\hat{l}(\mu) \geq 7.823$, then \[l(\gamma) < \dfrac{2\pi}{\hat{l}(\mu)^2 - 28.78}.\] 
\end{theorem}

The following result quantifies Thurston's Dehn filling theorem and guarantees that if a slope $\sigma \subset \partial M$ is sufficiently long, then the hyperbolic metric on $M$ deforms to a complete hyperbolic metric on the surgered manifold $M(\sigma)$, in which the core $v$ of the surgery solid torus becomes a simple geodesic.

\begin{theorem}
\label{theorem:fpscore}
    Let $M$ be a hyperbolic 3-manifold with torus boundary components and let $\sigma$ be a slope along a boundary component. 
    If $\hat{l}(\sigma)>10.1$, then $M(\sigma)$ is hyperbolic and the core $v$ of the surgery solid torus is a simple geodesic in $M(\sigma)$.
\end{theorem}

\begin{proof}
    The statement follows from the proof of \cite[Theorem~7.28]{fps}.
    Notice that the tools (\cite[Theorem~5.17]{fps}) used in this proof when filling all the cusps of a hyperbolic $3$-manifold still apply in our setting, irrespective of the fact that only a single cusp is being filled. 
\end{proof}

The meridian of $v$ is precisely the slope $\sigma$. 
Therefore we can apply Theorem \ref{theorem:fps} to relate the normalised length $\hat{l}(\sigma)$ of the filling slope $\sigma$ and the hyperbolic length $l(v)$ of the core curve $v$ of the filling.
\section{Determining a value for \texorpdfstring{$\mathcal{C}(K)$}{C(K)}}
\label{sec:realising}

We now turn to the proof of Theorem \ref{theorem:C(K)}. 
Let $K \subset \mathbb{S}^3$ be a prime knot whose exterior $\sphere_K$ is not a graph manifold. 
Since hyperbolic knots have already been studied \cite{wakelin}, we may also assume that $K$ is a satellite knot. 

Suppose that there is an orientation-preserving homeomorphism $f: \mathbb{S}^3_K(p/q) \to \mathbb{S}^3_{K'}(p/q)$ for some other knot $K' \subset \mathbb{S}^3$, where $|q| > 2$. 
By Proposition \ref{prop:JSJsurgery}, $K'$ is also a satellite knot and there are JSJ pieces $Y$ and $Y'$ of $\sphere_{K}$ and $\sphere_{K'}$ such that $Y(p/qt^2)$ and $Y'(p/q{t'}^2)$ are the surgered pieces of $\sphere_K(p/q)$ and $\sphere_{K'}(p/q)$, respectively (where $t,t' \geq 1$). 

The main part of our argument lies in proving Proposition \ref{proposition:swapping}. 
This tells us that if \mbox{$|q| > \max\{R(K), S(K)\}$}, then the initial orientation-preserving homeomorphism \mbox{$f: \mathbb{S}^3_K(p/q) \to \mathbb{S}^3_{K'}(p/q)$} restricts to one between the surgered pieces, so that \mbox{$f (Y(p/qt^2)) = Y'(p/q{t'}^2)$}. 

To prove Proposition \ref{proposition:swapping}, we will go through each possible case from Theorem \ref{theorem:budney} for $Y'$. 
Firstly, note that $Y'$ cannot be the exterior of a torus knot, otherwise $\sphere_{K'}(p/q) \cong Y'(p/qt'^2)$ would contain no JSJ tori. 
Secondly, if $Y'$ is a composing space, then $\sphere_{K'}(p/q)$ contains the filling $\sphere_{P'}(p/qt'^2)$ of a composing pattern $P'$, which must be homeomorphic to a knot exterior if $f(Y(p/qt^2)) \neq Y'(p/qt'^2)$; according to \cite[Lemma 5.4]{sorya}, this cannot happen for $|q|>1$. 
This leaves us with two cases: either $Y'$ is hyperbolic or $Y'$ is a cable space.

\subsection{The hyperbolic case}

Suppose now that $Y'$ is hyperbolic.

Recall that for a hyperbolic JSJ piece $X \subset \sphere_K$, we defined
\[\mathfrak{s}(L_X) =  \Bigg\lfloor\sqrt{6\sqrt{3}\bigg(\frac{2\pi}{\sys(\sphere_{L_X})}+28.78\bigg)}\Bigg\rfloor,\] 
where $L_X$ is the unique link corresponding to $X$ in the satellite construction of $K$.
We claim that taking $|q|>\max\{32,\mathfrak{s}(L_X)\}$ obstructs the filling $Y'(p/qt'^2)$ of our hyperbolic JSJ piece $Y'$ from becoming homeomorphic to another hyperbolic JSJ piece $X \subset \mathbb{S}^3_K \setminus Y$, i.e. $f(X) \neq Y'(p/q{t'}^2)$. 

Note that the definition of $\mathfrak{s}(L)$ for a hyperbolic link $L$ only depends on the homeomorphism type of $\sphere_L$, so we may also define $\mathfrak{s}(X)$ for any hyperbolic 3-manifold $X$, with $\sys(X)$ taking the place of $\sys(\sphere_L)$, to obtain a more general statement. 

\begin{lemma} 
\label{lemma:hyperbolicbound}
    Let $Y'$ be a hyperbolic JSJ piece of a knot exterior. 
    Consider a non-trivial slope $\sigma' = a'/b'$ on the outermost boundary component of $Y'$. 
    Let $X$ be any hyperbolic 3-manifold.
    
    If $|b'| > \max \left\{32, \mathfrak{s}(X) \right\}$, then $Y'(a'/b') \not\cong X$. 
    \begin{proof}
        By the 6-theorem \cite{agol,lackenby_6}, the meridional slope $\mu=1/0$ of the outermost boundary component of $Y'$ has length $l(\mu)\leq 6$. 
        Applying Lemma \ref{lemma:length_inequalities} to the slopes $\mu=1/0$ and $\sigma'=a'/b'$ gives
        $$\hat{l}(\sigma') \geq \frac{|b'|}{\hat{l}(\mu)} \geq \frac{|b'|}{\sqrt{6\sqrt{3}}}.$$
        Suppose that $|b'| > \max\{32,\mathfrak{s}(X)\}$. 
        Then
        $$\hat{l}(\sigma') \geq \frac{33}{\sqrt{6\sqrt{3}}} \geq 10.1 \geq 7.823$$
        and hence we can apply Theorem \ref{theorem:fpscore} and Theorem \ref{theorem:fps}. 
        We also have 
        $$\hat{l}(\sigma') \geq \frac{\mathfrak{s}(X)+1}{\sqrt{6\sqrt{3}}} > \sqrt{\frac{2\pi}{\sys(X)}+28.78}$$
        and so Theorem \ref{theorem:fps} gives us the following bound on the length of the core curve $v'$ of the filling $Y'(a'/b')$: 
        $$l(v') < \frac{2\pi}{\hat{l}(\sigma')^2-28.78} < \sys(X).$$ 
        Since $v'$ is a geodesic in $Y'(a'/b')$ with length strictly less than $\sys(X)$, we have $Y'(a'/b') \not \cong X$.
    \end{proof} 
\end{lemma} 

We now apply Lemma \ref{lemma:hyperbolicbound} to the setting of Proposition \ref{proposition:swapping}. 
Recall that
\[S(K) = \max\{32, \mathfrak{s}(L_X) \ | \ X \in \mathcal{X}, \ \partial \sphere_K \not \subset X\},\] 
where $\mathcal{X}$ is the set of hyperbolic JSJ pieces $X \subset \sphere_K$. 

\begin{proof}[Proof of Proposition \ref{proposition:swapping} (hyperbolic case)]
    Suppose that $f(X) = Y'(p/qt'^2)$ for some JSJ piece $X$ of $\sphere_K(p/q)$ which is not its surgered piece. 
    Since $|qt'^2|\geq |q| > 2$, $X$ must have been a hyperbolic JSJ piece of $\sphere_K$, so $X \in \mathcal{X}$ and $\partial \sphere_K \not \subset X$. 
    Let $L_X$ be the unique link corresponding to $X$ in the satellite construction of $K$. 
    By Lemma \ref{lemma:hyperbolicbound}, we have $|q| \leq |qt'^2| \leq \mathfrak{s}(L_X) \leq S(K)$. 
\end{proof}

\subsection{The cable case}
\label{subsection:cablecase}

We are left with the case when $Y'$ is a cable space. 

Let $X$ be a hyperbolic JSJ piece of $\sphere_K$ and let $L_X=L_0 \cup L_1 \cup \ldots \cup L_{m-1}$ be the unique link corresponding to $X$, where $L_0$ is the outermost component as usual. 
Recall that we defined
\[\mathfrak{r}(L_X) = \Bigg\lfloor\sqrt{3}\max\{0,l(\mu_i) \;|\; L_i \subset L, i\neq 0\}\Bigg\rfloor.\]

If $Y'$ is a cable space and the homeomorphism \mbox{$f: \mathbb{S}^3_K(p/q) \to \mathbb{S}^3_{K'}(p/q)$} does not restrict to one between the surgered pieces, then $Y'(p/qt'^2)$ is an irreducible Seifert fibre space and thus its preimage is the exterior $\sphere_J \subset \sphere_K$ of a torus knot $J$ by Theorem \ref{theorem:budney}. 
Let $P$ be the pattern such that $K=P(J)$. 
The JSJ torus $\partial \sphere_{J} \subset \sphere_{K}(p/q)$ is mapped by $f$ to the boundary of $Y'(p/qt'^2)$ in $\sphere_{K'}(p/q)$, which also bounds the exterior $\sphere_{J'} \subset \sphere_{K'}$ of a knot $J'$. 
It follows that $K'$ can be described as a satellite $P'(J')$ such that the homeomorphism $f$ restricts to $f(\sphere_P(p/q))=\sphere_{J'}$ and $f(\sphere_J)=\sphere_{P'}(p/q)$. 

We claim that taking $|q| > R(K)$ obstructs this situation from happening, thus completing the proof of Proposition \ref{proposition:swapping}. 
A key ingredient is our assumption that the exterior of $K$ is not a graph manifold: in particular, since $\sphere_J$ is a torus knot exterior, the pattern space $\sphere_P$ must contain a hyperbolic JSJ piece $X$. 
We'll  perform a clever filling of $\mathbb{S}^3_P$ that realises \cite[Proposition 5.8]{sorya} in our context. 

\begin{lemma}
\label{lemma:notSFS}
    Let $K$ be a satellite knot by a pattern $P=Q \cup U$. 
    Let $\mathcal{X}_P$ be the set of hyperbolic JSJ piece of $\sphere_P$ and, for each $X \in \mathcal{X}_P$, let $L_X$ denote the unique link corresponding to $X$ in the satellite construction of $K$. 
    Let $a/b$ be a slope on the boundary component $T$ of $\mathbb{S}^3_P$ corresponding to $U$. 
    
    If $\mathcal{X}_P \neq \varnothing$ and $|b|> \max\{1,\mathfrak{r}(L_{X}) \ | \ X \in \mathcal{X}_P\}$, then the filling $\sphere_P(T; a/b)$ is not a Seifert fibre space.
\end{lemma}

\begin{proof}
    Let $L_X = L_0 \cup \ldots \cup L_{m-1}$ be the link corresponding to an $X \in \mathcal{X}_P$, where $L_0$ is outermost.

    If $\sphere_P(T; a/b)$ is a Seifert fibre space, then no $X \in \mathcal{X}_P$ can be a JSJ piece of $\sphere_P(T; a/b)$. 
    Since $|b|>1$, the proof of \cite[Proposition 5.8]{sorya} implies that there is some $X \in \mathcal{X}_P$ and some $i \neq 0$ such that $\partial \nu (L_i)$ either is $T$ itself or is glued to the exterior of an iterated cable pattern containing $T$ inside $\sphere_P$. 
    In both cases, we see that $\sphere_P(T; a/b)$ contains a manifold $X(L_i; a/bc^2), X \in \mathcal{X}_P$, such that $|bc^2|\geq|b|>\mathfrak{r}(L_{X})$.
    
    We claim that the assumption $|b| > \mathfrak{r}(L_{X})$ ensures that the filling $X(L_i;a/bc^2)$ remains hyperbolic. 
    Applying Lemma \ref{lemma:length_inequalities} to the distinct pair of slopes $\sigma=a/bc^2$ and $\mu_i = 1/0$ along $L_i$ gives
    \[l(\sigma) \geq \frac{2\sqrt{3}|bc^2|}{l(\mu_i)} \geq \frac{2\sqrt{3}(\mathfrak{r}(L_X)+1)}{l(\mu_i)} > \frac{2\sqrt{3}(\sqrt{3}l(\mu_i))}{l(\mu_i)} = 6.\]
    By the 6-theorem, $X(L_i;a/bc^2)$ is hyperbolic. Therefore $\sphere_P(T; a/b)$ cannot be a Seifert fibre space.
\end{proof}

We'll apply this to solve our last remaining case of Proposition \ref{proposition:swapping}. 
Recall that 
$$R(K) = \max\{1, \mathfrak{r}(L_X) \ | \ X \in \mathcal{X} \},$$ 
where $\mathcal{X}$ is the set of hyperbolic JSJ pieces $X \subset \sphere_K$. 

\begin{proof}[Proof of Proposition \ref{proposition:swapping} (cable case)]
    Write $K=P(J)$ and $K'=P'(J')$ as in the discussion at the beginning of this subsection. 
    Let $T = \sphere_P \cap \sphere_J$. 
    Since $Y'$ is a cable space, $P'$ has winding number $w'\neq0$. 
    By \cite[Lemma 4.3]{sorya}, the preimage of the meridian $\mu_{J'}$ of $J'$ has slope $y/q{w'}^2$ along $T$, for some integer $y$, in the coordinates given by the link component $U$ of $P = Q \cup U$ which corresponds to $T$. 
    Let $M = \sphere_P(T; y/q{w'}^2)$. 
    Observe that 
            \begin{align*}
                M(p/q) &\cong \mathbb{S}^3_{J'}(\mu_{J'}) = \mathbb{S}^3, \\
                M(1/0) &\cong \mathbb{S}^3_U(\mu_{J'}) = L(y,qw'^2). 
            \end{align*} 
    By Lemma \ref{lemma:notSFS}, the manifold $M$ is not a Seifert fibre space. 
    However, $M(p/q)$ and $M(1/0)$ are cyclic surgeries for $M$ such that $|q| > 1$, contradicting the cyclic surgery theorem \cite{cgls}.
\end{proof}

\subsection{Proof of Theorem \ref{theorem:C(K)}} \label{subsection:proofmaintheorem}

Having completed the proof of Proposition \ref{proposition:swapping}, we see that the orientation-preserving homeomorphism $f:\sphere_K(p/q) \to \sphere_{K'}(p/q)$ restricts to \mbox{$f (Y(p/qt^2)) = Y'(p/q{t'}^2)$}. 
We claim that in fact $Y=Y'$ as JSJ pieces. 

\begin{proposition}
\label{proposition:Y=Y'}
    Let $K$ be a prime satellite knot whose exterior is not a graph manifold. 
    Suppose that there is an orientation-preserving homeomorphism $f: \mathbb{S}^3_K(p/q) \to \mathbb{S}^3_{K'}(p/q)$ for some knot $K'$.
    
    Let $Y$ and $Y'$ be the JSJ pieces such that $f$ restricts to a slope-preserving homeomorphism between the surgered pieces $Y(p/qt^2)$ and $Y'(p/q{t'}^2)$. 
    Let $L$ and $L'$ be the links corresponding to $Y$ and $Y'$ in the satellite constructions of $K$ and $K'$, respectively. 
    
    If $|q| > Q(K)$, then $Y=Y'$, in the sense that $Y \cong Y'$ and $L=L'$. 
    \begin{proof}
        We run through the cases in Theorem \ref{theorem:budney}. 
        Note that $Y$ cannot be a torus knot exterior, as the exterior of $K$ is assumed to contain at least one hyperbolic JSJ piece. 
        If $Y$ is a composing space or a cable space, then \cite[Sections 6.2 and 6.4]{sorya} imply that $Y = Y'$. 
        If $Y$ is a hyperbolic JSJ piece, then taking $|q|> Q(K) = \max\{34,\mathfrak{q}(L)\}$ and applying \cite[Proposition~4.8]{wakelin} implies that $Y = Y'$.
    \end{proof}
\end{proposition} 

Finally, we return to the caveat that the outermost JSJ pieces of $\sphere_K$ and $\sphere_{K'}$ may or may not be cable spaces that become solid tori after filling, as in case (ii) of Proposition \ref{prop:JSJsurgery}. 
This can be resolved by considering cosmetic surgeries. 

\begin{proof}[Proof of Theorem \ref{theorem:C(K)}] 
    Let $K$ be a prime satellite knot whose exterior is not a graph manifold. 
    Suppose that there is an orientation-preserving homeomorphism $f: \mathbb{S}^3_K(p/q) \rightarrow \mathbb{S}^3_{K'}(p/q)$ for some knot $K'$. 
    By Proposition \ref{proposition:swapping}, this restricts to a slope-preserving homeomorphism between the surgered pieces: $f(Y(p/qt^2)) = Y'(p/q{t'}^2)$, where $t,t'\geq1$. 
    By Proposition \ref{proposition:Y=Y'}, it follows that $Y = Y'$. 
    Therefore we can write $K$ and $K'$ as (possibly trivial) cables of the same non-trivial knot $\widehat{K}$ whose exterior has outermost piece $Y=Y'$. 
    Observe that we can write: 
    \begin{align*}     
        \mathbb{S}^3_K(p/q) &\cong \mathbb{S}^3_{\widehat{K}}(p/qt^2), \\ 
        \mathbb{S}^3_{K'}(p/q) &\cong \mathbb{S}^3_{\widehat{K}}(p/q{t'}^2). 
    \end{align*}
    If $t \neq t'$, then $\widehat{K}$ has a pair of distinct cosmetic surgery slopes of the same sign, contradicting \cite[Theorem~1.2]{nw}. 
    Therefore $t=t'$ and, by \cite[Lemma~6.2]{sorya}, we deduce that $K=K'$. 
\end{proof}
\section{Winding number zero} 
\label{sec:windingzero}

The goal of this section is to prove Theorem \ref{theorem:C(K)usingt}. 
Recall that $K$ will now be a satellite knot for which every satellite description $K=P(J)$ is by a pattern $P$ with winding number zero. 
We will show that a different strategy can be used to obstruct the swapping of JSJ pieces in an orientation-preserving homeomorphism \mbox{$f: \mathbb{S}^3_K(p/q) \to \mathbb{S}^3_{K'}(p/q)$}: instead of assuming that $|q|>\max\{R(K),S(K)\}$ as in Proposition \ref{proposition:swapping}, we will take $|q|>\max\{2,T(K)\}$ and prove Proposition \ref{proposition:swappingusingt}. 
The new bound $T(K)$ arises from showing that for such a knot $K$, the swapping of JSJ pieces can only occur for finitely many possibilities, which can sometimes be identified and avoided directly.

\subsection{Splicifiable knots}

Recall that a \emph{nullhomologous Rolfsen $t$-twist} on a knot $J$ refers to performing $-1/t$-surgery along an unknot in $\mathbb{S}^3$ which is nullhomologous in the exterior of $J$, thus adding $t$ full twists to $J$ in this location. 
We call the integer $t$ involved in this process the \emph{nullhomologous Rolfsen twist coefficient}. 

\begin{definition}
\label{definition:splicifiable}
    Let $K=P(J)$ be a satellite knot. 
    If the pattern $P=Q\cup U$ has winding number zero, $Q$ is unknotted and $J$ can be unknotted by some nullhomologous Rolfsen $t$-twist, then we say that $K$ is \emph{$t$-splicifiable} with respect to the pair $(P,J)$. 
\end{definition} 

The reason for this choice of terminology is illustrated in the following result. 

\begin{proposition}
\label{proposition:splicifiable}
    Let $K=P(J)$ be a satellite knot by a pattern $P$ with winding number zero. 

    There exists a satellite knot $K' = P'(J')$ and an orientation-preserving homeomorphism \mbox{$f:\mathbb{S}^3_K(p/q) \to \mathbb{S}^3_{K'}(p'/q')$} with $f(\mathbb{S}^3_P(p/q)) = \mathbb{S}^3_{J'}$ and $f(\mathbb{S}^3_J) = \mathbb{S}^3_{P'}(p'/q')$ if and only if \mbox{$|p|=|p'|=1$} and $K$ is $\pm q'$-splicifiable with respect to $(P,J)$. 
    
    Furthermore, $\mathbb{S}^3_K(\mp 1/q) \cong \mathbb{S}^3_{K'}(\mp 1/q')$ is homeomorphic to the splice of the knot exteriors $\sphere_J$ and $\sphere_{J'}$ and $K'$ is $\pm q$-splicifiable with respect to $(P',J')$. 
\end{proposition} 

\begin{remark}
    Brakes \cite{brakes} provides a general method for constructing an orientation-preserving homeomorphism \mbox{$f: \mathbb{S}^3_{P(J)}(p/q) \rightarrow \mathbb{S}^3_{{P'(J')}}(p'/q')$} such that $f(\mathbb{S}^3_P(p/q)) = \mathbb{S}^3_{J'}$ and $f(\mathbb{S}^3_J)=\mathbb{S}^3_{P'}(p'/q')$. 
    Proposition \ref{proposition:splicifiable} says that if such a manifold is obtained from a knot $K=P(J)$, where $P$ is a pattern with winding number zero, then Brakes' construction is in fact the only possibility. 
    We will soon see that this restricts the existence of non-characterising slopes for $K$ with large denominator. 
\end{remark}

The following lemma is key to the proof of Proposition \ref{proposition:splicifiable}. 

\begin{lemma} 
\label{lemma:splice}
    Let $K=P(J)$ and $K'=P'(J')$ be satellite knots and let $p/q$ and $p'/q'$ be non-trivial slopes along $K$ and $K'$, respectively. 
    Suppose that there is an orientation-preserving homeomorphism $f: \mathbb{S}^3_K(p/q) \to \mathbb{S}^3_{K'}(p'/q')$ with \mbox{$f(\mathbb{S}^3_P(p/q))=\mathbb{S}^3_{J'}$} and \mbox{$f(\mathbb{S}^3_J)=\mathbb{S}^3_{P'}(p'/q')$}. 

    Then $P$ has winding number $w=0$ if and only if $P'$ has winding number $w'=0$. 
    Moreover, if these winding numbers are zero, then $|p|=|p'|=1$, $P$ and $P'$ are both unknotted patterns and $\mathbb{S}^3_K(p/q) \cong \mathbb{S}^3_{K'}(p'/q')$ is homeomorphic to the splice of the knot exteriors $\mathbb{S}^3_J$ and $\mathbb{S}^3_{J'}$. 
    \begin{proof}
        Write $P=Q\cup U$ and $P'=Q'\cup U'$, where $Q$ and $Q'$ are knots inside the respective solid tori $\mathbb{S}^3_U$ and $\mathbb{S}^3_{U'}$. 
        From our satellite constructions, we have the following identifications between meridians and longitudes: 
        \[K=P(J) \iff
            \begin{cases}
                \mu_J = \lambda_U \\ 
                \lambda_J = \mu_U
            \end{cases}
        \hspace{-1em};\]
        \[K'=P'(J') \iff 
            \begin{cases}
                \mu_{J'} = \lambda_{U'} \\ 
                \lambda_{J'} = \mu_{U'}
            \end{cases} 
        \hspace{-1em}.\]
        Given the homeomorphism $f$ between fillings, recall from \cite[Lemma~4.3]{sorya} that the meridians and longitudes of companion knots can be expressed as follows (for some integers $y$ and $y'$): 
        \[f(\mathbb{S}^3_P(p/q))=\mathbb{S}^3_{J'} \implies 
            \begin{cases}
                \mu_{J'}=y'\mu_U + q'{w'}^2\lambda_U \\ 
                \lambda_{J'}=q{w'}^2\mu_U + p\lambda_U 
            \end{cases}
        \hspace{-1em};\]
        \[f(\mathbb{S}^3_{J})=\mathbb{S}^3_{P'}(p'/q') \implies 
            \begin{cases}
                \mu_{J}=y\mu_{U'} + qw^2\lambda_{U'} \\ 
                \lambda_J=q'w^2\mu_{U'} + p'\lambda_{U'} 
            \end{cases}
        \hspace{-1em}.\]
        We observe that 
        \[w'=0 \iff 
            \begin{cases}
                \mu_{J'}=\mu_U \\ 
                \lambda_{J'}=\lambda_U
            \end{cases}
            \iff 
            \begin{cases}
                \mu_{J'}=\lambda_J \\ 
                \lambda_{J'}=\mu_J 
            \end{cases}
            \iff 
            \begin{cases}
                \lambda_{U'}=\lambda_J \\
                \mu_{U'}=\mu_J 
            \end{cases}
        \iff w=0. 
        \]
        Moreover, this gluing map produces precisely the splice of the knot exteriors $\mathbb{S}^3_J$ and $\mathbb{S}^3_{J'}$. 

        Consider the manifolds $\mathbb{S}^3_P(U;\mu_{J'})$ and $\mathbb{S}^3_{P'}(U';\mu_J)$. 
        Since $\mu_{J'}=\mu_U$ and $\mu_J=\mu_{U'}$, these manifolds are in fact the knot exteriors $\mathbb{S}^3_Q$ and $\mathbb{S}^3_{Q'}$, respectively. 
        These both have non-trivial $\mathbb{S}^3$ fillings: 
        \begin{align*}
            \mathbb{S}^3_Q(p/q) \cong \sphere_{J'}(\mu_{J'}) &\cong \mathbb{S}^3, \\ 
            \mathbb{S}^3_{Q'}(p'/q') \cong \sphere_J(\mu_J) &\cong \mathbb{S}^3.
        \end{align*}
        We conclude that $P$ and $P'$ must be unknotted patterns and that $|p|=|p'|=1$, as required. 
    \end{proof}
\end{lemma} 

Combining this lemma with Brakes' construction, we complete the proof of Proposition \ref{proposition:splicifiable}. 

\begin{proof}[Proof of Proposition \ref{proposition:splicifiable}]
    Let $K=P(J)$ be $\pm q'$-splicifiable with respect to $(P,J)$ and write \mbox{$P=Q \cup U$}. 
    Since $Q$ is unknotted and $P$ has winding number zero, we can now perform a nullhomologous Rolfsen $\pm q$-twist to $U$ along $Q$, for any choice of coefficient with $|q|\geq1$, to obtain a knot $J'$ with \mbox{$f(\mathbb{S}^3_P(\mp1/q)) = \mathbb{S}^3_{J'}$}. 
    Since $J$ can be unknotted by a Rolfsen $\pm q'$-twist, there must also be a pattern $P'$ for which $f(\mathbb{S}^3_J) = \mathbb{S}^3_{P'}(\mp1/q')$. 
    Thus Brakes' construction gives a $\pm q$-splicifiable knot $K' = P'(J')$ as in the statement of the theorem. 
    Lemma \ref{lemma:splice} gives the other direction. 
\end{proof}

By specialising $\pm q$ to be a nullhomologous Rolfsen twist coefficient which unknots the companion of a splicifiable knot, we may obtain non-characterising slopes. 

\begin{corollary}
\label{corollary:nonchar1/t}
    Let $K=P(J)$ be a satellite knot which is $-t$-splicifiable with respect to $(P,J)$. 
    
    Then there exists a satellite knot $K'=P'(J')$ which is  $-t$-splicifiable with respect to $(P',J')$ such that there is an orientation-preserving homeomorphism $f: \mathbb{S}^3_K(1/t) \to \mathbb{S}^3_{K'}(1/t)$ with $f(\mathbb{S}^3_P(1/t)) = \mathbb{S}^3_{J'}$ and $f(\mathbb{S}^3_J) = \mathbb{S}^3_{P'}(1/t)$.

    Furthermore, if $K \neq K'$, then $1/t$ is a non-characterising slope for $K$ and $K'$. 
\end{corollary}

This generalises \cite[Theorem 1.8]{wakelin}, which corresponds to the special case when $K=P(J)$ is a multiclasped Whitehead double of a double twist knot, demonstrating that this process can be used to realise non-characterising slopes with arbitrarily high denominator. 

An extra assumption on $K$ ensures that this construction always gives non-characterising slopes.

\begin{lemma}
\label{lemma:nonchar}
    Let $K=P(J)$ be a satellite knot which is $-t$-splicifiable with respect to $(P,J)$. Let $J'$ be the knot such that $\mathbb{S}^3_P(1/t) \cong \sphere_{J'}$.

    If $K$ cannot be expressed as a satellite of $J'$, then $1/t$ is a non-characterising slope for $K$.
\end{lemma}

\begin{proof}
    By Corollary \ref{corollary:nonchar1/t}, there is a knot $K'$ such that there is a homeomorphism $f: \sphere_K(1/t) \to \sphere_{K'}(1/t)$ with $f(\mathbb{S}^3_P(1/t)) = \mathbb{S}^3_{J'} \subset \mathbb{S}^3_{K'}$. However, by our assumption, $\mathbb{S}^3_{J'} \not \subset \mathbb{S}^3_{K}$; hence $\mathbb{S}^3_K$ and $\mathbb{S}^3_{K'}$ have different JSJ decompositions. Therefore $K' \neq K$ and, by Corollary \ref{corollary:nonchar1/t}, $1/t$ is a non-characterising slope.
\end{proof}

\subsection{Maximal nullhomologous Rolfsen twist coefficient}

We are now ready to define the constant $T(K)$ appearing in Theorem \ref{theorem:C(K)usingt}. 

First, we observe that the construction leading to the non-characterising slopes in Corollary \ref{corollary:nonchar1/t} bounds an unknotting nullhomologous Rolfsen twist coefficient. To the best of the authors' knowledge, this is the first such proof which does not require extra conditions on either the knot or the twist. 

\begin{corollary}
\label{corollary:maxtwist}
    Let $J$ be a non-trivial knot which can be unknotted by a nullhomologous Rolfsen $t$-twist for some $t \in \mathbb{Z}$. There there is a maximal possible value $\mathfrak{t}(J)>0$ for $|t|$. 
    \begin{proof}
        Set $K=W^n(J)$ to be a multiclasped Whitehead double of $J$, where $|n|$ is larger than any parameter $|m|, |t|$ appearing in a double twist knot exterior $\mathbb{S}^3_{T^m_t} \subset \mathbb{S}^3_J$. 
        Observe that $W^n$ is an unknotted pattern and that the surgered piece is homeomorphic to a double twist knot exterior, $\mathbb{S}^3_{W^n}(\mp1/q) \cong \mathbb{S}^3_{T^n_{\pm q}}$. By construction, $K$ cannot be expressed as a satellite of $T^n_{\pm q}$. 

        Suppose for contradiction that no such $\mathfrak{t}(J)$ exists. Then there are infinitely many $q \in \Z$ such that $J$ can be unknotted by a single nullhomologous Rolfsen $\pm q$-twist. 
        Thus $K$ is $\pm q$-splicifiable with respect to $(W^n,J)$ for infinitely many $q \in \Z$. 
        By Lemma \ref{lemma:nonchar}, we obtain infinitely many non-characterising slopes $\mp 1/q$ for $K$. 
        This contradicts \cite[Theorem 1.1]{sorya}. 
    \end{proof}
\end{corollary} 

Given such a knot $J$, we call the integer $\mathfrak{t}(J)$ its \emph{maximal nullhomologous Rolfsen twist coefficient}. 

\begin{definition}
\label{def:T(K)}
    Define $T(K) := \max\{0,\mathfrak{t}(J) \;|\; K \text{ is splicifiable with respect to a pair }(P,J)\}$.
\end{definition} 

We are now in a position to prove Theorem \ref{theorem:C(K)usingt}.

\subsection{Proof of Theorem \ref{theorem:C(K)usingt}} 

We begin with the proof of Proposition \ref{proposition:swappingusingt}, which uses the bound $T(K)$ to obstruct the swapping of JSJ pieces in a homeomorphism between surgeries. Recall that we are in the case where $K$ is a satellite knot for which every satellite description $K=P(J)$ is by a pattern $P$ with winding number zero. 

\begin{proof}[Proof of Proposition \ref{proposition:swappingusingt}]
    Suppose that there is an orientation-preserving homeomorphism \mbox{$f: \mathbb{S}^3_K(p/q) \to \mathbb{S}^3_{K'}(p/q)$} for some slope $p/q$ with $|q|> \max\{2,T(K)\}$. Since we have $|q|>2$, Proposition \ref{prop:JSJsurgery} ensures that the JSJ pieces of the surgered manifolds are well-defined. 
    
    Suppose that the surgered pieces of $\mathbb{S}^3_K(p/q)$ and $\mathbb{S}^3_{K'}(p/q)$ are not mapped to one another by $f$. Then, by the same reasoning as in Subsection \ref{subsection:cablecase}, we can describe the knots $K, K'$ as satellites $P(J), P'(J')$, respectively, such that $f(\mathbb{S}^3_P(p/q))= \mathbb{S}^3_{J'}$ and $f(\mathbb{S}^3_J)=\mathbb{S}^3_{P'}(p/q)$. By Proposition \ref{proposition:splicifiable}, the knot $K$ must be $\pm q$-splicifiable, which implies that $J$ can be unknotted by a nullhomologous Rolfsen $\pm q$-twist. By Definition \ref{def:T(K)}, we have $|q| \leq \mathfrak{t}(J) \leq T(K)$, a contradiction.
\end{proof}

The orientation-preserving homeomorphism $f$ thus restricts to one between the surgered pieces. It remains to deduce that $K=K'$. 

\begin{proof}[Proof of Theorem \ref{theorem:C(K)usingt}]
    Let $K$ be a satellite knot such that for every choice of satellite description $K=P(J)$, the pattern $P$ has winding number zero. 
    Suppose that there is an orientation-preserving homeomorphism \mbox{$f: \sphere_K(p/q) \to \sphere_{K'}(p/q)$} for some knot $K'$ and slope $p/q$ with $|q|>\max\{Q(K), T(K)\}$. 
    By Proposition \ref{proposition:swappingusingt}, this restricts to a slope-preserving homeomorphism between the surgered pieces. 
    Since $K$ is a satellite knot of hyperbolic type, we can now apply Proposition \ref{proposition:Y=Y'} and the cosmetic surgery argument used in the proof of Theorem \ref{theorem:C(K)} to deduce that $K=K'$. 
\end{proof}  
\section{Examples} 
\label{sec:applications}

We conclude this article with a series of examples exhibiting the utility of our main results. 

We will begin by showcasing Theorem \ref{theorem:C(K)} through some illustrative examples. Recall that the bound \mbox{$\mathcal{C}(K)=\max\{Q(K), R(K), S(K)\}$} only depends on the hyperbolic JSJ pieces of the knot exterior $\mathbb{S}^3_K$. We will find a value for this $\mathcal{C}(K)$ using the computer programme SnapPy \cite{cdgw_SnapPy}.

For certain knots, Theorem \ref{theorem:C(K)usingt} will give a refinement $\mathcal{C}(K)=\max\{Q(K),T(K)\}$. This depends on the maximal nullhomologous Rolfsen twist coefficient of a companion for $K$, which is generally harder to compute, but we will use the fact that $T(K)$ is known (often to be just $0$ or $1$) in many cases.

\subsection{Examples of Theorem \ref{theorem:C(K)}} 
We will first give an example which simply demonstrates how to compute the bound in Theorem \ref{theorem:C(K)}. We will then see how this is affected by making modifications to the knot.

\begin{example}
\label{example:borromean} 
    Let $K = B(W(3_1), 4_1 \# 6_1)$ be the satellite knot of hyperbolic type constructed by splicing the Borromean rings $B$ with $W(3_1)$ (the Whitehead double of the right-handed trefoil) and $4_1 \# 6_1$ (the connected sum of the figure-eight knot and the stevedore knot). 

    By Theorem \ref{theorem:C(K)}, only the hyperbolic JSJ pieces of $\mathbb{S}^3_K$ contribute to $\mathcal{C}(K)=\max\{Q(K),R(K),S(K)\}$. We have
    \begin{align*}
        Q(K) &= \max\{34, \mathfrak{q}(B)\} = \max\{34,18\} = 34; \\ 
        R(K) &= \max\{1, \mathfrak{r}(B), \mathfrak{r}(W), \mathfrak{r}({4_1}), \mathfrak{r}({6_1})\} = \max\{1,2,2,0,0\} = 2; \\ 
        S(K) &= \max\{32, \mathfrak{s}(W), \mathfrak{s}({4_1}), \mathfrak{s}({6_1})\} = \max\{32, 18, 18, 22\} = 32; 
    \end{align*}
    which gives $\mathcal{C}(K)=\max\{34,2,32\}=34$. 
\end{example} 

\begin{remark}
\label{remark:cable}
    Let $\widehat{K}$ be a knot of hyperbolic type with hyperbolic outermost JSJ piece $Y$ and corresponding link $L_Y$. Let $K=C_{r,s}(\widehat{K})$ be any cable of $\widehat{K}$. The only change from the bound $\mathcal{C}(\widehat{K})$ to the bound $\mathcal{C}(K)$ is the extra contribution of $\mathfrak{s}(L_Y)$, as $Y$ is no longer the outermost JSJ piece. However, it is easy to see from the formulae that $\mathfrak{s}(L_Y)\leq \mathfrak{q}(L_Y)$. Therefore we can in fact take $\mathcal{C}(K)=\mathcal{C}(\widehat{K})$. 
\end{remark} 

\begin{example}
\label{example:cable}
    Let $K = C_{1,2}(B(W(3_1), 4_1 \# 6_1))$ be the $(1,2)$-cable of the knot in Example \ref{example:borromean}. Then we can take $\mathcal{C}(K)=34$ by Remark \ref{remark:cable}. 
\end{example} 

\begin{remark}
\label{remark:twisting}
    Let $L = L_0 \cup U^{m-1}$ be a hyperbolic link and consider any link obtained by adding a nullhomologous Rolfsen twist to $L$ along a component of $U^{m-1}$. Performing such a twist does not change the homeomorphism type of the link exterior, so its systole is unchanged and hence both $\mathfrak{q}(L)$ and $\mathfrak{s}(L)$ are unaffected. However, $\mathfrak{r}(L)$ is defined in terms of the meridians of the link components, so any such change to $L$ may affect this. 
\end{remark}

\begin{example}
\label{example:twisting} 
    Let $K = B_{-5,2}(W_{-7}(3_1), 4_1 \# 6_1)$ be the knot constructed in almost the same way as the one in Example \ref{example:borromean}, but with a $-7$-twisted Whitehead link and with $B_{-5,2}$ denoting the Borromean rings twisted along two of its unlink components $-5$ and $2$ times, respectively. 
    Then we can compute 
        \[R(K) = \max\{1, \mathfrak{r}(B_{-5,2}), \mathfrak{r}(W_{-7}), \mathfrak{r}({4_1}), \mathfrak{r}({6_1})\} = \max\{1, 24, 36, 0, 0\} = 36\]
    and take $\mathcal{C}(K)=\max\{Q(K), R(K), S(K)\}=\max\{34,36,32\}=36$. 
\end{example}

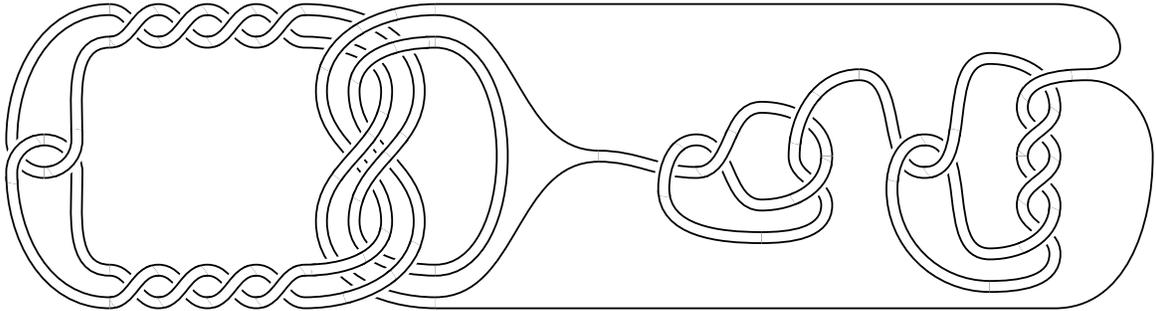
\begin{figure}[htbp!]
    \centering 
    \resizebox{\textwidth}{!}{\centering
\begin{tikzpicture} [squarednode/.style={rectangle, draw=black, fill=white, thick, minimum size=40pt}]
    \begin{knot}[end tolerance = 1pt]  
        \flipcrossings{2,4,5,7,8,10,11,13,15,18,16,19,21,23,25,27,29,31}
    
        
        \strand[black, double=white, thick, double distance=4pt] (-1.5,-1.75) 
        to [out=left, in=down] (-2.75,-1) 
        to [out=up, in=down] (-1.75,1) 
        to [out=up, in=right] (-3.5,2.25) 
        to [out=left, in=right] (-4.25,1.75)
        to [out=left, in=right] (-5,2.25)
        to [out=left, in=right] (-5.75,1.75) 
        to [out=left, in=right] (-6.5,2.25); 
        \strand[black, double=white, thick, double distance=4pt] (-6.5,2.25) 
        to [out=left, in=up] (-8,0.25)
        to [out=down, in=left] (-7.5,-0.25); 
        \strand[black, double=white, thick, double distance=4pt] (-7.5,-0.25) 
        to [out=right, in=down] (-7,0.25)
        to [out=up, in=left] (-6.5,1.75); 
        \strand[black, double=white, thick, double distance=4pt] (-6.5,1.75) 
        to [out=right, in=left] (-5.75,2.25) 
        to [out=right, in=left] (-5,1.75)
        to [out=right, in=left] (-4.25,2.25) 
        to [out=right, in=left] (-3.5,1.75) 
        to [out=right, in=up] (-2.25,1) 
        to [out=down, in=up] (-3.25,-1) 
        to [out=down, in=left] (-1.5,-2.25); 
        
        \strand[black, double=white, thick, double distance=4pt] (-1.5,1.75) 
        to [out=left, in=up] (-2.75,1) 
        to [out=down, in=up] (-1.75,-1) 
        to [out=down, in=right] (-3.5,-2.25)
        to [out=left, in=right] (-4.25,-1.75) 
        to [out=left, in=right] (-5,-2.25) 
        to [out=left, in=right] (-5.75,-1.75)
        to [out=left, in=right] (-6.5,-2.25); 
        \strand[black, double=white, thick, double distance=4pt] (-6.5,-2.25) 
        to [out=left, in=down] (-8,-0.25)
        to [out=up, in=left] (-7.5,0.25); 
        \strand[black, double=white, thick, double distance=4pt] (-7.5,0.25) 
        to [out=right, in=up] (-7,-0.25)
        to [out=down, in=left] (-6.5,-1.75); 
        \strand[black, double=white, thick, double distance=4pt] (-6.5,-1.75) 
        to [out=right, in=left] (-5.75,-2.25) 
        to [out=right, in=left] (-5,-1.75) 
        to [out=right, in=left] (-4.25,-2.25) 
        to [out=right, in=left] (-3.5,-1.75) 
        to [out=right, in=down] (-2.25,-1) 
        to [out=up, in=down] (-3.25,1) 
        to [out=up, in=left] (-1.5,2.25); 

        \strand[black, double=white, thick, double distance=4pt] (-1.5,-1.75) 
        to [out=right, in=right] (-1.5,1.75); 


        \strand[black, thick] (-1.5,2.165) 
        to [out=right, in=left] (1,0.085); 

        \strand[black, thick] (-1.5,-2.165) 
        to [out=right, in=left] (1,-0.085); 

        \strand[black, thick] (-1.5,2.335) 
        to [out=right, in=left] (8,2.335) 
        to [out=right, in=up] (9,1.665) 
        to [out=down, in=right] (8.5,1.335); 

        \strand[black, thick] (-1.5,-2.335) 
        to [out=right, in=left] (8,-2.335) 
        to [out=right, in=down] (9.5,0) 
        to [out=up, in=right] (8.5,1.165); 
        

        \strand[black, double=white, thick, double distance=4pt] (1,0) 
        to [out=right, in=left] (2.5,-0.25) 
        to [out=right, in=left] (3.5,0.75)
        to [out=right, in=up] (4.5,0); 
        
        \strand[black, double=white, thick, double distance=4pt] (4.5,0) 
        to [out=down, in=right] (3.5,-0.75)
        to [out=left, in=right] (2.5,0.25)
        to [out=left, in=up] (2,-0.5) 
        to [out=down, in=left] (3.5,-1.25); 
        
        \strand[black, double=white, thick, double distance=4pt] (3.5,-1.25) 
        to [out=right, in=down] (4.5,-0.75)
        to [out=up, in=down] (4,0)
        to [out=up, in=left] (5,1.25); 


        \strand[black, double=white, thick, double distance=4pt] (5,1.25) 
        to [out=right, in=left] (6,-0.25) 
        to [out=right, in=left] (7,1.5)
        to [out=right, in=up] (8,0.75) 
        to [out=down, in=up] (7.5,0); 
        
        \strand[black, double=white, thick, double distance=4pt] (7.5,0) 
        to [out=down, in=up] (8,-0.75)
        to [out=down, in=right] (7,-1.5)
        to [out=left, in=right] (6,0.25)
        to [out=left, in=up] (5.5,-0.5) 
        to [out=down, in=left] (7,-2); 
        
        \strand[black, double=white, thick, double distance=4pt] (7,-2) 
        to [out=right, in=down] (8,-1.5)
        to [out=up, in=down] (7.5,-0.75)
        to [out=up, in=down] (8,0)
        to [out=up, in=down] (7.5,0.75) 
        to [out=up, in=left] (8.5,1.25); 
    \end{knot}
\end{tikzpicture} } 
    \caption{The knot $B_{-5,2}(W_{-7}(3_1), 4_1 \# 6_1)$ from Example \ref{example:twisting}.}
    \label{figure:borromean}
\end{figure}

\subsection{Examples of Theorem \ref{theorem:C(K)usingt}} 

In many cases, the alternative swapping obstruction used in the proof of Theorem \ref{theorem:C(K)usingt} allows us to refine the realisation of $\mathcal{C}(K)$. Not only do the non-characterising slopes in Corollary \ref{corollary:nonchar1/t} give a lower bound on the optimal value for $\mathcal{C}(K)$, but Theorem \ref{theorem:C(K)usingt} may also improve the realisation obtained by Theorem \ref{theorem:C(K)}. 

Although it might not always be easy to check whether a knot can be unknotted by a single nullhomologous Rolfsen twist, nor to determine the constant $T(K)$, we will show that $T(K)=0$ for certain satellites of knots with large signature and $T(K)\leq1$ for certain satellites of knots which are composite or fibred. Furthermore, we'll see that when $T(K) \geq Q(K)$, we obtain an optimal value for $\mathcal{C}(K)$.

\subsubsection{Satellites of knots with large signature} 

In the following situation, we will see that $T(K)=0$.

\begin{corollary}
\label{corollary:T=0}
    Let $K$ be a satellite knot such that for every choice of satellite description $K=P(J)$, the winding number of $P$ is zero but $K$ is not splicifiable with respect to $(P,J)$.
    
    If $|q|> Q(K)$, then $p/q$ is a characterising slope for $K$. 
    \begin{proof}
        By definition, we have $T(K)=0$ because there is no satellite description $K=P(J)$ such that $K$ is splicifiable with respect to $(P,J)$. By Theorem \ref{theorem:C(K)usingt}, we have that every slope $p/q$ with $|q|> \max\{Q(K), T(K)\}=Q(K)$ is characterising for $K$. 
    \end{proof}
\end{corollary}

Recall that the \emph{surgery description number} $sd(K)$ of a knot $K$ is defined to be the minimum number of regions required to unknot $K$ via nullhomologous Rolfsen twists. 
Note that $\mathfrak{t}(K)$ is only defined when $sd(K)=1$. 
The surgery description number is related to several other knot invariants \cite{aikrt}. Here, we observe that the signature of a knot gives a lower bound for its surgery description number. 

\begin{lemma}
\label{lemma:signature}
    Let $K$ be a knot. Then its signature $\sigma(K)$ satisfies
    $\frac{|\sigma(K)|}{2} \leq sd(K)$.
    \begin{proof}
         The signature $\sigma(K)$ is a lower bound for twice the topological 4-genus, $2 g^{\text{TOP}}_4(K)$ \cite{kauffman-taylor}. Let $g_{a}(K)$ be the minimal difference between the genera of a Seifert surface $F$ for $K$ and a subsurface $F' \subset F$ bounded by a knot $K'$ with Alexander polynomial $\Delta_{K'}=1$. We have $g^{\text{TOP}}_4(K) \leq g_{a}(K)$ since a knot $K'$ with $\Delta_{K'}=1$ has $g^{\text{TOP}}_4(K)=0$ \cite{freedman}. Performing the nullhomologous Rolfsen twists relating $K$ to the unknot in succession, we obtain a sequence of knots $K_i, i=0,\ldots, sd(K)$, where $K_0 = K$ and $K_{sd(K)}=U$ is the unknot. Using \cite[Theorem 1.1]{mccoy_galg}, we see that
         \[g_a(K)=|g_a(K_0)-g_a(K_{sd(K)})| = \sum_{i=0}^{sd(K)-1} |g_a(K_i)-g_a(K_{i+1})| \leq sd(K).\]
         Combining this with the earlier inequalities gives the result.
    \end{proof}
\end{lemma}

By definition, the companion $J$ of a satellite knot $P(J)$ that is splicifiable with respect to $(P,J)$ must have $sd(J)=1$. 
Combining Corollary \ref{corollary:T=0} with Lemma \ref{lemma:signature} yields the following. 

\begin{corollary}
\label{corollary:signature}
    Let $K$ be a satellite knot such that for every choice of satellite description $K=P(J)$, the winding number of $P$ is zero and the companion $J$ has signature satisfying $|\sigma(J)|\geq 4$.

    If $|q|> Q(K)$, then $p/q$ is a characterising slope for $K$. 
    \begin{proof}
        Write $K=P(J)$. Since $|\sigma(J)|\geq4$, Lemma \ref{lemma:signature} implies that $J$ cannot be unknotted via a nullhomologous Rolfsen twist so $K$ is not splicifiable with respect to $(P,J)$. This being true for every satellite description $K=P(J)$, we apply Corollary \ref{corollary:T=0}. 
    \end{proof}
\end{corollary}

Below is an example of a knot for which this result yields a better bound than the one given by Theorem \ref{theorem:C(K)}. 

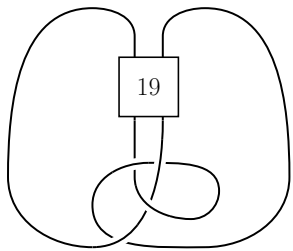
\begin{figure}
    \centering
    \resizebox{0.25\textwidth}{!}{\begin{tikzpicture} [squarednode/.style={rectangle, draw=black, fill=white, thick, minimum size=40pt}]
\begin{scope}
    \begin{knot}[end tolerance=1pt] 
    \flipcrossings{1,3,4} 
        \strand[white, double=black, thick, double distance=1pt] (-1,-1.25) 
        to [out=left, in=down] (-2.5,0) 
        to [out=up, in=left] (-1,3)
        to [out=right, in=up] (-0.25,2.5)
        to [out=down, in=up] (-0.25,1); 
        \strand[white, double=black, thick, double distance=1pt] (-0.25,1) 
        to [out=down, in=up] (-0.25,0)
        to [out=down, in=left] (0.751,-0.75)
        to [out=right, in=down] (1.25,-0.25) 
        to [out=up, in=right] (0,0.25); 
        \strand[white, double=black, thick, double distance=1pt] (0,0.25)
        to [out=left, in=up] (-1,-0.5)
        to [out=down, in=left] (1,-1.25); 
        \strand[white, double=black, thick, double distance=1pt] (0.25,01) 
        to [out=up, in=down] (0.25,1) 
        to [out=up, in=down] (0.25,2.5) 
        to [out=up, in=left] (1,3) 
        to [out=right, in=up] (2.5,0) 
        to [out=down, in=right] (1,-1.25); 
        \strand[white, double=black, thick, double distance=1pt] (0.25,01) to [out=down, in=right] (-1,-1.25); 
    \end{knot} 
    \node[squarednode, fill=white, opacity=1, draw=black, scale=0.75] at (0,1.6) {\huge $19$}; 
\end{scope} 
\end{tikzpicture} } 
    \caption{A knot with ``small'' systole and ``large'' signature.} 
    \label{figure:knotsignature}
\end{figure}

\begin{example}
    Let $J$ be the hyperbolic knot pictured in Figure \ref{figure:knotsignature} and let $K=W(J)$ be its Whitehead double. 
    We have $Q(K) = \max\{34,\mathfrak{q}(W)\}=34$ and $R(K)=\max\{1,\mathfrak{r}(W)\}=1$. 
    Moreover, SnapPy tells us that $\sys(\sphere_J) \approx 0.0141687$, so $S(K)=\max\{32,\mathfrak{s}(J)\}=70$.
    Theorem \ref{theorem:C(K)} then gives a realisation of $\mathcal{C}(K)$ as $\max\{34,1,70\}=70$. However, we also have $\sigma(J)=-38$, so $T(K)=0$ and our bound can be improved to $\max\{Q(K),T(K)\}=\max\{34,0\}=34$ by Corollary \ref{corollary:signature}. 
\end{example}

\subsubsection{Satellites of composite and fibred knots} 

Lackenby \cite{lackenby_sd} showed that for any composite or fibred knot $J$, if $sd(J)=1$ then in fact $\mathfrak{t}(J)=1$. This yields the following corollary for certain satellites of such knots. 

\begin{corollary}
\label{corollary:compositefibred}
    Let $K$ be a satellite knot such that for every choice of satellite description $K=P(J)$, the winding number of $P$ is zero, and if $P$ is an unknotted pattern, then the companion $J$ is either composite or fibred.

    If $|q|> Q(K)$, then $p/q$ is a characterising slope for $K$. 
    \begin{proof}
        By assumption, $K$ may be splicifiable with respect to a pair $(P,J)$ only when $J$ is either a composite knot or a fibred knot. Hence any knot $J$ that contributes to $T(K)$ has maximal nullhomologous Rolfsen twist coefficient $\mathfrak{t}(J)=1$ according to \cite{lackenby_sd}. Therefore $Q(K)\geq 34 >1 \geq T(K)$ and we apply Theorem \ref{theorem:C(K)usingt}.
    \end{proof} 
\end{corollary}

Corollary \ref{corollary:compositefibred} provides more examples of knots for which the bound obtained from this result is an improvement on the bound coming from Theorem \ref{theorem:C(K)}. First, we make the following simple observation. 

\begin{lemma}
\label{lemma:q<s}
Let $K$ be a satellite knot of hyperbolic type. If 
    \[\sys(X)\leq\frac{12\sqrt{3} \pi }{Q(K)^2-172.68 \sqrt{3}}\] for every hyperbolic JSJ piece $X$ which is not outermost in $\sphere_K$, then $Q(K) \leq S(K)$. 
    \begin{proof}
        The hypothesis implies that $Q(K) \leq \mathfrak{s}(L_X)$ for every non-outermost hyperbolic JSJ piece $X$, where $L_X$ is the link corresponding to $X$ in the satellite construction of $K$. Hence $Q(K) \leq S(K)$. 
    \end{proof}
\end{lemma}

This allows us to construct examples where the bound from Theorem \ref{theorem:C(K)} is $\mathcal{C}(K)=S(K)$ but Theorem \ref{theorem:C(K)usingt} gives an improved bound $\mathcal{C}(K)=Q(K)$. 

\begin{example}
\label{example:fibred} 
    Let $J$ be a fibred hyperbolic knot and let $K=W(J)$ be its Whitehead double. Applying Corollary \ref{corollary:compositefibred}, we obtain a realisation of $\mathcal{C}(K)$ as $Q(K) = 34$. By Lemma \ref{lemma:q<s}, this is an improvement of the realisation obtained by Theorem \ref{theorem:C(K)} whenever
    \[\sys(\sphere_J) \leq \frac{12\sqrt{3} \pi }{34^2-172.68 \sqrt{3}} \approx 0.0762003.\]

    For instance, take $J$ to be the fibred pretzel knot $P(-2,-77,77)$ \cite{gabai_pretzel}. We have \[\sys(\sphere_J) \approx 0.0035737 \leq 0.0762003.\] Whilst Corollary \ref{corollary:compositefibred} realises the bound $\mathcal{C}(K)$ as $Q(K)=34$, Theorem \ref{theorem:C(K)} realises $\mathcal{C}(K)$ as $\max\{Q(K), R(K), S(K)\}= \max\{34, 1, 136\} = 136 > 34$. 
\end{example}

\begin{example}
\label{example:composite}
    Let $J$ be the connected sum of two simple knots $J_1$ and $J_2$ and let $K=W(J)$ be its Whitehead double. There are three possible satellite descriptions of $K$, each corresponding to a JSJ torus of $\sphere_K$: $W(J), P_2(J_1)$ and $P_1(J_2)$, where $P_i$ is the composing pattern $W(J_i) \cup U$ for $i=1,2$. Since the winding number of $W$ is zero, the winding number of $P_i$ is also zero. Furthermore, the component $W(J_i)$ of $P_i$ is knotted. Therefore $K$ satisfies the conditions of Corollary \ref{corollary:compositefibred}, and we obtain a realisation of $\mathcal{C}(K)$ as $34$. If $J_1$ or $J_2$ is hyperbolic, suppose that the condition of Lemma \ref{lemma:q<s} is satisfied. Then, as in Example \ref{example:fibred}, this is an improvement of the bound coming from Theorem \ref{theorem:C(K)}. 
 \end{example}

\subsubsection{Optimal bounds} 

The bound from Theorem \ref{theorem:C(K)} is unlikely to be optimal due to the nature of its construction. However, in some cases the refined bound from Theorem \ref{theorem:C(K)usingt} is truly optimal. 

\begin{corollary}
\label{corollary:optimal}
    Let $K$ be a satellite knot such that for every choice of satellite description $K=P(J)$, $K$ is splicifiable with respect to $(P,J)$. Suppose that $T(K) \geq Q(K)$.

    If $|q| > T(K)$, then $p/q$ is a characterising slope for $K$ and this is the optimal such bound. 
    \begin{proof}
        First, observe that $1/T(K)$ is a non-characterising slope for $K$. Lemma \ref{lemma:splice} tells us that any non-characterising slope with larger denominator would have to correspond to an orientation-preserving homeomorphism between surgeries which restricts to one between the surgered pieces. However, if $|q|> Q(K)$, then no such non-characterising slope can exist. Hence our bound is optimal when $T(K) \geq Q(K)$. 
    \end{proof} 
\end{corollary} 

Whilst $Q(K)$ is computable \cite{hw}, it is generally harder to find an explicit value for $T(K)$. In all of our previous examples, we had $T(K)<Q(K)$. Suppose that $K=P(J)$ is a satellite of a simple knot $J$ by a hyperbolic pattern $P$ such that $K$ is splicifiable with respect to $(P,J)$. If it is known that $\mathfrak{t}(J)>1$, then one can follow Lackenby's algorithm in \cite{lackenby_algorithm} to find the exact value of $\mathfrak{t}(J)$ and hence $T(K)$. The following example shows that this can be made arbitrarily high, so that $T(K) \geq Q(K)$. 

\begin{example} 
    Let $K=W^n(T^m_t)$ be a multiclasped Whitehead double of a double twist knot with $\max\{m,t\}>1$. Both $1/m$ and $1/t$ are non-characterising slopes for $K$ which can be realised by nullhomologous Rolfsen twists unknotting $T^m_t$ \cite[Theorem~1.8]{wakelin}. Since $\mathfrak{t}(T^m_t) \geq \max\{m,t\}>1$, we may choose $m,n,t$ such that $T(K) \geq Q(K) = \max\{34, \mathfrak{q}(W^n)\}$ and we can apply Corollary \ref{corollary:optimal}. For instance, if $n=1$, then for any $m,t \geq 34$, we have that $T(K)\geq Q(K)=34$. Thus $T(K)$, which can be obtained by Lackenby's algorithm, is optimal.
\end{example}  

\vspace{1em}
\bibliographystyle{amsalpha}
\bibliography{bib} 

\end{document}